# MAXIMAL DIVISORIAL SETS IN ARC SPACES

SHIHOKO ISHII

1 2

ABSTRACT. In this paper we introduce a maximal divisorial set in the arc space of a variety. The generalized Nash problem is reduced to a translation problem of the inclusion of two maximal divisorial sets. We study this problem and show a counter example to the most natural expectation even for a non-singular variety.

## 1. INTRODUCTION

In [8], Nash posed a problem: if the set of the families of arcs through the singularities on a variety (these families are called the Nash components) corresponds bijectively to the set of essential divisors of resolutions of the singularities. This problem is affirmatively answered for some 2-dimensional singularities by A. Reguera, M. Lejeune-Jalabert, C. Plénat and P. Popescu-Pampu [7], [9], [10], [11], [12] and toric singularities of arbitrary dimension by S. Ishii and J. Kollár [6]. For non-normal toric variety of arbitrary dimension the answer is also affirmative ([5]). On the other hand this problem is negatively answered in general. The paper [6] gives a counter example of dimension greater than or equal to 4. Therefore the Nash problem should be changed to the problem to determine the divisors corresponding to the Nash components.

We can generalize this problem into the characterization problem for valuations corresponding to the irreducible components of contact loci $\text{Cont}^{\geq m}(\mathfrak{a})$ which are introduced by L. Ein, R. Lazarsfeld and M. Mustaţă ([2]). In this paper we introduce the maximal divisorial set $C_X(v)$ in the arc space of a variety $X$ corresponding to a divisorial valuation $v$. A fat irreducible component of a contact locus is a maximal divisorial set. In order to characterize the valuations corresponding to the irreducible components of a contact locus, it is essential to translate the inclusion relation between two maximal divisorial sets to a relation

[1]partially supported by Grant-In-Aid of Ministry of Science and Education in Japan
[2]Mathematics Subject Classification 2000: Primary 14J17, Secondary 14M99

between the corresponding divisorial valuations. The most natural candidate as this relation is the value-inequality relation, i.e., $v(f) \leq v'(f)$ for every regular function $f$ on the affine variety $X$. Actually, if the variety $X$ and the valuations $v, v'$ are toric, we have the equivalence: $v(f) \leq v'(f)$ for every regular function $f$ on $X \Leftrightarrow C_X(v) \supset C_X(v')$. For non-toric valuations, we show that this equivalence does not hold even on $\mathbb{C}^2$.

This paper is organized as follows: In the second section we introduce the maximal divisorial set corresponding to a divisorial valuation and show the basic properties. In the third section we show some basic properties of a contact locus and formulate a generalized Nash problem. In the forth section we show an example of divisorial valuations $v, v'$ over $\mathbb{C}^2$ such that the value-inequality relation does not imply the inclusion of the maximal divisorial sets corresponding to $v$ and $v'$.

The author is grateful to the members of Singularity Seminar at Nihon University for useful suggestions and encouragement.

In this paper, a variety is always an irreducible reduced separated scheme of finite type over $\mathbb{C}$.

## 2. Maximal divisorial sets in the arc space

**Definition 2.1.** Let $X$ be a scheme of finite type over $\mathbb{C}$ and $K \supset \mathbb{C}$ a field extension. A morphism $\alpha : \operatorname{Spec} K[[t]] \longrightarrow X$ is called an *arc* of $X$. We denote the closed point of $\operatorname{Spec} K[[t]]$ by 0 and the generic point by $\eta$.

For $m \in \mathbb{N}$, a morphism $\beta : \operatorname{Spec} K[t]/(t^{m+1}) \longrightarrow X$ is called an *m-jet* of $X$. Denote the space of arcs of $X$ by $X_\infty$ and the space of $m$-jets of $X$ by $X_m$.

The arc space $X_\infty$ and $m$-jet space $X_m$ are characterized by the following property:

**Proposition 2.2.** *Let $X$ be a scheme of finite type over $\mathbb{C}$. Then, for an arbitrary $\mathbb{C}$-scheme $Y$,*

$$\operatorname{Hom}_\mathbb{C}(Y, X_\infty) \simeq \operatorname{Hom}_\mathbb{C}(Y \widehat{\times}_{\operatorname{Spec} \mathbb{C}} \operatorname{Spec} \mathbb{C}[[t]], X)$$

*where $Y \widehat{\times}_{\operatorname{Spec} \mathbb{C}} \operatorname{Spec} \mathbb{C}[[t]]$ means the formal completion of $Y \times_{\operatorname{Spec} \mathbb{C}} \operatorname{Spec} \mathbb{C}[[t]]$ along the subscheme $Y \times_{\operatorname{Spec} \mathbb{C}} \{0\}$, and for $m \in \mathbb{N}$*

$$\operatorname{Hom}_\mathbb{C}(Y, X_m) \simeq \operatorname{Hom}_\mathbb{C}(Y \times_{\operatorname{Spec} \mathbb{C}} \operatorname{Spec} K[t]/(t^{m+1}), X).$$

**2.3.** By thinking of the case $Y = \operatorname{Spec} K$ for an extension field $K$ of $\mathbb{C}$, we see that $K$-valued points of $X_\infty$ correspond to arcs $\alpha : \operatorname{Spec} K[[t]] \longrightarrow X$ bijectively. Based on this, we denote the $K$-valued point corresponding to an arc $\alpha : \operatorname{Spec} K[[t]] \longrightarrow X$ by the same symbol $\alpha$. The canonical projection $X_\infty \longrightarrow X$, $\alpha \mapsto \alpha(0)$ is denoted by $\pi_X$. If there is no risk of confusion, we write just $\pi$.

A morphism $\varphi : Y \longrightarrow X$ of varieties induces a canonical morphism $\varphi_\infty : Y_\infty \longrightarrow X_\infty$, $\alpha \mapsto \varphi \circ \alpha$.

The concept "thin" in the following is first introduced in [2] and a "fat arc" is introduced and studied in [5].

**Definition 2.4.** Let $X$ be a variety over $\mathbb{C}$. We say an arc $\alpha : \operatorname{Spec} K[[t]] \longrightarrow X$ is *thin* if $\alpha$ factors through a proper closed subset of $X$. An arc which is not thin is called a *fat arc*.

An irreducible closed subset $C$ in $X_\infty$ is called a *thin set* if the generic point of $C$ is thin. An irreducible closed subset in $X_\infty$ which is not thin is called a *fat set*.

**Definition 2.5.** Let $\alpha : \operatorname{Spec} K[[t]] \longrightarrow X$ be a fat arc of a variety $X$ and $\alpha^* : \mathcal{O}_{X,\alpha(0)} \longrightarrow K[[t]]$ the local homomorphism induced from $\alpha$. By Proposition 2.5, (i) in [5], $\alpha^*$ is extended to the injective homomorphism of fields $\alpha^* : K(X) \longrightarrow K((t))$, where $K(X)$ is the rational function field of $X$. Define a function $v_\alpha : K(X) \setminus \{0\} \longrightarrow \mathbb{Z}$ by

$$v_\alpha(f) = \operatorname{ord}_t \alpha^*(f).$$

Then, $v_\alpha$ is a valuation of $K(X)$. We call it the *valuation corresponding to* $\alpha$.

**Definition 2.6.** Let $X$ be a variety, $g : X_1 \longrightarrow X$ a proper birational morphism from a normal variety $X_1$ and $E \subset X_1$ an irreducible divisor. Let $f : X_2 \longrightarrow X$ be another proper birational morphism from a normal variety $X_2$. The birational map $f^{-1} \circ g : X_1 \dashrightarrow X_2$ is defined on a (nonempty) open subset $E^0$ of $E$. The closure of $(f^{-1} \circ g)(E^0)$ is well defined. It is called the *center* of $E$ on $X_2$.

We say that $E$ appears in $f$ (or in $X_2$), if the center of $E$ on $X_2$ is also a divisor. In this case the birational map $f^{-1} \circ g : X_1 \dashrightarrow X_2$ is a local isomorphism at the generic point of $E$ and we denote the birational transform of $E$ on $X_2$ again by $E$. For our purposes $E \subset X_1$ is identified with $E \subset X_2$. (Strictly speaking, we should be talking about the corresponding *divisorial valuation* instead.) Such an equivalence class is called a *divisor over* $X$.

**Definition 2.7.** A valuation $v$ on the rational function field $K(X)$ of a variety $X$ is called a *divisorial valuation over* $X$ if $v = q \operatorname{val}_E$ for some $q \in \mathbb{N}$ and a divisor $E$ over $X$. The center of a divisor $E$ is called *the center of* the valuation $v = q \operatorname{val}_E$. A fat arc $\alpha$ of $X$ is called a *divisorial arc* if $v_\alpha$ is a divisorial valuation over $X$. A fat set is called a *divisorial set* if the generic point is a divisorial arc.

**Definition 2.8.** For a divisorial valuation $v$ over a variety $X$, define the maximal divisorial set corresponding to $v$ as follows:

$$C_X(v) := \overline{\{\alpha \in X_\infty \mid v_\alpha = v\}},$$

where $\overline{\{\ \}}$ is the Zariski closure in $X_\infty$.

**Proposition 2.9.** *Let $\varphi : Y \longrightarrow X$ be a proper birational morphism of varieties and $v$ a divisorial valuation over $X$. Then,*

  (i) *$v$ is also a divisorial valuation over $Y$ and $C_X(v) = \overline{\varphi_\infty(C_Y(v))}$.*
  (ii) *If an open subset $U \subset Y$ intersects the center of the valuation $v$ on $Y$, Then $C_X(v) = \overline{\varphi_\infty(C_U(v))}$.*

*Proof.* The first assertion of (i) follows from constructing a suitable proper birational morphism dominating $Y$ on which the corresponding divisor appears. The inclusion $\supset$ is obvious. For the opposite inclusion, take a fat arc $\alpha \in X_\infty$ such that $v_\alpha = v$. As the image $\alpha(\eta)$ is the generic point of $X$, it is in the open subset on which $\varphi$ is isomorphic. Hence, $\alpha(\eta)$ is lifted on $Y$. Then, by the valuative criterion of properness, $\alpha$ is uniquely lifted to $\tilde\alpha \in Y_\infty$. We obtain that $\alpha = \varphi(\tilde\alpha)$ and $v_{\tilde\alpha} = v$, which imply $\alpha \in \varphi_\infty(C_Y(v))$. For the statement (ii), let $\alpha : \operatorname{Spec} K[[t]] \longrightarrow Y$ be a fat arc of $Y$ such that $v_\alpha = v$. Then, $\alpha(\eta)$ is the generic point of $Y$ and $\alpha(0)$ is the generic point of the center of $v$ on $Y$. One can see that both are on $U$, therefore $\alpha \in U_\infty$. □

In the following section, we will see that a maximal divisorial set is irreducible.

**Example 2.10.** For toric varieties, we use the terminologies in [3]. Let $X$ be an affine toric variety defined by a cone $\sigma$ in $N \simeq \mathbb{Z}^n$. Let $M$ be the dual of $N$ and $\langle\ ,\ \rangle : N \times M \longrightarrow \mathbb{Z}$ the canonical pairing. A $\mathbb{C}$-algebra $\mathbb{C}[M]$ consists of linear combinations over $\mathbb{C}$ of monomials $x^u$ for $u \in M$. A point $v \in \sigma \cap N$ gives a toric valuation $v$ by $v(x^u) = \langle v, u \rangle$. An irreducible locally closed subset $T_\infty(v)$ is defined in [4] as follows:

$$T_\infty(v) = \{\alpha \in X_\infty \mid \alpha(\eta) \in T,\ \operatorname{ord}_t \alpha^*(x^u) = \langle v, u \rangle \text{ for } u \in M\},$$

where $T$ denotes the open orbit and also the torus acting on $X$ and $\alpha^* : \mathbb{C}[\sigma^\vee \cap M] \longrightarrow K[[t]]$ is the ring homomorphism corresponding to $\alpha$. Then, $C_X(v) = \overline{T_\infty(v)}$. This is proved as follows: First, take the generic point $\alpha \in C_X(v)$. Then, $\operatorname{ord}_t \alpha^*(x^u) = v(x^u) = \langle v, u \rangle$. Therefore, $\alpha \in T_\infty(v)$, which yield $C_X(v) \subset \overline{T_\infty(v)}$. Conversely, let $\beta \in \overline{T_\infty(v)}$ be the generic point. Then, by the upper semi-continuity (see for example [5, Proposition 2.7]) $v_\beta(f) \leq v_\alpha(f) = v(f)$ for every $f \in \mathbb{C}[\sigma^\vee \cap M]$. But, since $v$ is toric, it is the minimal valuation satisfying $v_\beta(x^u) = v(x^u)$ for every $u \in \sigma^\vee \cap M$. Therefore, $v_\beta = v$, which implies $\beta \in C_X(v)$.

## 3. Contact loci

**Definition 3.1.** Let $\psi_m : X_\infty \longrightarrow X_m$ be the canonical projection to the space of $m$-jets $X_m$. A subset $C \subset X_\infty$ is called a *cylinder* if there is a constructible set $S \subset X_m$ for some $m \in \mathbb{N}$ such that
$$C = \psi_m^{-1}(S).$$

For an arc $\alpha : \operatorname{Spec} K[[t]] \longrightarrow X$ of an affine variety $X = \operatorname{Spec} A$, we always denote by $\alpha^*$ the ring homomorphism $A \longrightarrow K[[t]]$ corresponding to $\alpha$.

Let $X$ be a variety and $\mathfrak{a}$ an ideal sheaf on $X$. For an arc $\alpha : \operatorname{Spec} K[[t]] \longrightarrow X$, there is an open affine subset $U \subset X$ such that $\alpha$ factors through $U$. We define $\operatorname{ord}_t \alpha^*(\mathfrak{a})$ as follows:
$$\operatorname{ord}_t \alpha^*(\mathfrak{a}) = \operatorname{ord}_t \alpha^*(\Gamma(U, \mathfrak{a})).$$

**Definition 3.2** ([2]). For an ideal sheaf $\mathfrak{a}$ on a variety $X$, we define
$$\operatorname{Cont}^m(\mathfrak{a}) = \{\alpha \in X_\infty \mid \operatorname{ord}_t \alpha^*(\mathfrak{a}) = m\}$$
and
$$\operatorname{Cont}^{\geq m}(\mathfrak{a}) = \{\alpha \in X_\infty \mid \operatorname{ord}_t \alpha^*(\mathfrak{a}) \geq m\}.$$
These subset are called *contact loci* of an ideal $\mathfrak{a}$. The subset $\operatorname{Cont}^{\geq m}(\mathfrak{a})$ is closed and $\operatorname{Cont}^m(\mathfrak{a})$ is locally closed. Both are cylinders.

**Definition 3.3** ([2]). For a simple normal crossing divisor $E = \bigcup_{i=1}^s E_i$ on a non-singular variety $X$, we introduce the multi-contact loci for a multi-index $\nu \in \mathbb{Z}_{\geq 0}^s$:
$$\operatorname{Cont}^\nu(E) = \{\alpha \in X_\infty \mid \operatorname{ord}_t \alpha^*(I_{E_i}) = \nu_i\},$$
where $I_{E_i}$ is the defining ideal of $E_i$. The multi-contact locus $\operatorname{Cont}^\nu(E)$ is irreducible if it is not empty.

**Proposition 3.4.** *Let $v = q \operatorname{val}_{E_0}$ be a divisorial valuation over a variety $X$. Let $\varphi : Y \longrightarrow X$ be a resolution of the singularities of $X$ such that the irreducible divisor $E_0$ appears on $Y$. Then,*
$$C_X(v) = \overline{\varphi_\infty(\operatorname{Cont}^q(E_0))}.$$
*In particular, $C_X(v)$ is irreducible.*

*Proof.* By Proposition 2.9, it is sufficient to prove that
$$C_Y(v) = \overline{\mathrm{Cont}^q(E_0)}.$$
Let $\alpha \in \mathrm{Cont}^q(E_0)$ be the generic point, then $\alpha(0) = e$ is the generic point of $E_0$. Therefore, we obtain a local homomorphism
$$\alpha^* : \mathcal{O}_{Y,e} \longrightarrow K[[t]]$$
such that $\mathrm{ord}_t \alpha^*(\tau) = q$, where $\tau$ is a generator of the maximal ideal of $\mathcal{O}_{Y,e}$. Then, for every $f \in \mathcal{O}_{Y,e}$ written as $f = u\tau^n$, where $u$ is a unit in $\mathcal{O}_{Y,e}$, it follows that $\mathrm{ord}_t \alpha^*(f) = qn = q\,\mathrm{val}_{E_0}(f) = v(f)$. This implies $\alpha \in C_Y(v)$.

Conversely, if $\beta \in C_Y(v)$ is an arc such that $v_\beta = v$, then $\mathrm{ord}_t \beta^*(\tau) = q$, which yields $\beta \in \mathrm{Cont}^q(E_0)$. □

The following is the characterization of ideals which have the same contact loci.

**Proposition 3.5.** *Let $\mathfrak{a}, \mathfrak{b} \subset A$ be ideals on an affine variety $X = \mathrm{Spec}\,A$. Then, the following are equivalent:*

(i) $\overline{\mathfrak{a}} = \overline{\mathfrak{b}}$, *where $\overline{\mathfrak{a}}, \overline{\mathfrak{b}}$ are the integral closures of $\mathfrak{a}$ and $\mathfrak{b}$, respectively;*

(ii) $\mathrm{Cont}^m(\mathfrak{a}) = \mathrm{Cont}^m(\mathfrak{b})$ *for every $m \in \mathbb{N}$;*

(iii) $\mathrm{Cont}^{\geq m}(\mathfrak{a}) = \mathrm{Cont}^{\geq m}(\mathfrak{b})$ *for every $m \in \mathbb{N}$.*

*Proof.* By $\mathrm{Cont}^m(\mathfrak{a}) = \mathrm{Cont}^{\geq m}(\mathfrak{a}) \setminus \mathrm{Cont}^{\geq m+1}(\mathfrak{a})$, the equivalence: $(ii) \Leftrightarrow (iii)$ is clear.

For $(i) \Rightarrow (iii)$, it is sufficient to prove that $\mathrm{Cont}^{\geq m}(\mathfrak{a}) = \mathrm{Cont}^{\geq m}(\mathfrak{a}')$ ($m \in \mathbb{N}$) for an ideal $\mathfrak{a}'$ integral over $\mathfrak{a}$. As $\mathfrak{a}' \supset \mathfrak{a}$, it follows $\mathrm{Cont}^{\geq m}(\mathfrak{a}') \subset \mathrm{Cont}^{\geq m}(\mathfrak{a})$. To prove the opposite inclusion, take an arc $\alpha \in \mathrm{Cont}^{\geq m}(\mathfrak{a})$. It is sufficient to prove that $\mathrm{ord}_t \alpha^*(x) \geq m$ for every $x \in \mathfrak{a}'$. Assume $\mathrm{ord}_t \alpha^*(x) = d < m$. Since $x$ is integral over $\mathfrak{a}$, there is a relation

$$(3.5.1) \qquad x^n + a_1 x^{n-1} + a_2 x^{n-2} + \ldots + a_n = 0$$

with $a_i \in \mathfrak{a}^i$. Here, noting that $\mathrm{ord}_t \alpha^*(a_i) \geq mi$, we obtain $\mathrm{ord}_t \alpha^*(a_i x^{n-i}) > nd$ for $i \geq 1$ and $\mathrm{ord}_t \alpha^*(x^n) = nd$. Therefore the order of the left hand side of (3.5.1) is $nd$, which is a contradiction to the equality.

For $(ii) \Rightarrow (i)$, take an element $x \in \mathfrak{b}$. Given a discrete valuation ring $R_v \supset A$, we obtain a ring homomorphism

$$A \subset R_v \longrightarrow \widehat{R_v} = K[[t]],$$

where $K$ is the residue field of $R_v$ and $t$ is a generator of the maximal ideal of $R_v$. Let $\alpha$ be the arc corresponding to this ring homomorphism $A \longrightarrow K[[t]]$. Let $m = \mathrm{ord}_t \alpha^*(\mathfrak{a})$. Then $\alpha \in \mathrm{Cont}^m(\mathfrak{a}) = \mathrm{Cont}^m(\mathfrak{b})$, and therefore $\mathrm{ord}_t \alpha^*(\mathfrak{b}) = m$. Hence, $v(x) = \mathrm{ord}_t \alpha^*(x) \geq m = v(\mathfrak{a} R_v)$.

This shows that $x$ is integral over $\mathfrak{a}$, i.e., $\mathfrak{b} \subset \overline{\mathfrak{a}}$ by the valuative characterization of integrity. Similarly, we obtain $\mathfrak{a} \subset \overline{\mathfrak{b}}$. Now we have $\overline{\mathfrak{a}} = \overline{\mathfrak{b}}$ as required. □

It is proved in [2, Corollary 2.6] that a fat component of a contact locus is a divisorial set in case $X$ is non-singular. The following gives more precise information also for singular $X$.

**Proposition 3.6** ([1]). *Let $X$ be an affine variety $\operatorname{Spec} A$ and $\mathfrak{a} \subset A$ a non-zero ideal. For $m \in \mathbb{N}$, a fat irreducible component of $\operatorname{Cont}^{\geq m}(\mathfrak{a})$ is a maximal divisorial set and the number of such components is finite.*

Here, we note that, for a divisorial valuation $v$ over an affine variety $X = \operatorname{Spec} A$ and an ideal $\mathfrak{a} \subset A$, $C_X(v) \subset \operatorname{Cont}^{\geq m}(\mathfrak{a})$ if and only if $v(\mathfrak{a}) \geq m$.

**Problem 3.7. A generalized Nash problem.** *(Embedded version of Nash's problem in [2]) Let $X$ be an affine variety $\operatorname{Spec} A$ and $\mathfrak{a} \subset A$ an ideal. Determine the set*

$$\mathcal{V}^m(\mathfrak{a}) := \{v \mid C_X(v) \text{ is an irreducible component of } \operatorname{Cont}^{\geq m}(\mathfrak{a})\}.$$

**Problem 3.8. The Nash problem.** *Let $I_{\operatorname{Sing} X} \subset A$ be the defining ideal of the singular locus $\operatorname{Sing} X$ of $X$. Determine the set $\mathcal{V}^1(I_{\operatorname{Sing} X})$.*

**3.9.** Nash posed his problem in a different way ([8], see also [6]), but his problem is translated into the above problem. As one sees, the Nash problem is a special case of the generalized Nash problem. Nash predicted that the set $\mathcal{V}^1(I_{\operatorname{Sing} X})$ coincides with the set of the valuations of essential divisors ([8]). But, it is not true for four or higher dimensional case ( [6]). So, Nash's prediction is still open for 2 and 3 dimensional cases.

**3.10.** The generalized Nash problem is to determine maximal $C_X(v)$'s contained in $\operatorname{Cont}^{\geq m}(\mathfrak{a})$. Therefore it is essential to determine the relation of valuations $v, v'$ when $C_X(v) \supset C_X(v')$. The most natural candidate for the relation of $v$ and $v'$ is that $v(f) \leq v'(f)$ for every $f \in A$, which is denoted by:

$$v|_A \leq v'|_A.$$

If $X$ is a toric variety and $v, v'$ are toric divisorial valuation, then we obtain in [4]:

$$v|_A \leq v'|_A \iff C_X(v) \supset C_X(v').$$

**Lemma 3.11.** *Let $v, v'$ be divisorial valuations over $X$.*
  (i) *If $C_X(v) \supset C_X(v')$, then $v|_A \leq v'|_A$.*
  (ii) *Assume that $v$ is a toric valuation, then the converse also holds.*

*Proof.* Let $\alpha$ and $\alpha'$ be the generic points of $C_X(v)$ and $C_X(v')$, respectively. Then, $C_X(v) \supset C_X(v')$ yields $\mathrm{ord}_t\, \alpha^*(f) \leq \mathrm{ord}_t\, \alpha'^*(f)$ for every $f \in A$ by the upper semi-continuity. This yields $v|_A \leq v'|_A$.

Assume that $v$ is toric. Let $v'_0$ be the toric valuation determined by $v'_0(x) = v'(x)$ for every monomial $x$. Then, by the definition of $T_\infty(v'_0)$, it follows $\alpha' \in T_\infty(v'_0)$. Hence, $C_X(v'_0) \supset C_X(v')$. On the other hand, $v$ and $v'_0$ are both toric, it follows that $v|_A \leq v'_0|_A$ and therefore $C_X(v) \supset C_X(v'_0)$ by [4, Proposition 4.8]. $\square$

So we expect that, at least for a simple variety like $\mathbb{C}^n$, the above equivalence holds for arbitrary divisorial valuations $v, v'$. In the next section we consider this problem and will give a negative answer.

## 4. An example over $\mathbb{C}^2$

**Proposition 4.1.** *For a divisorial valuations $v, v'$ over $X = \mathrm{Spec}\, A$, the following hold:*

(i) $C_X(v) = \overline{\bigcap_{f \in A \setminus \{0\}} \mathrm{Cont}^{v(f)}(f)}$;

(ii) $C_X(v)$ *is an irreducible component of* $\bigcap_{f \in A \setminus \{0\}} \mathrm{Cont}^{\geq v(f)}(f)$;

(iii) *the relation $v|_A \leq v'|_A$ holds if and only if the inclusion $C_X(v') \subset \bigcap_{f \in A \setminus \{0\}} \mathrm{Cont}^{\geq v(f)}(f)$ holds.*

*Proof.* For the proof of (i), take the generic point $\alpha$ of $C_X(v)$, then for every $f \in A \setminus \{0\}$ it follows that $\mathrm{ord}_t\, \alpha^*(f) = v(f)$ which means $\alpha \in \mathrm{Cont}^{v(f)}(f)$. On the other hand, let $\beta \in \bigcap_{f \in A \setminus \{0\}} \mathrm{Cont}^{v(f)}(f)$, then it follows $v_\beta = v$, which yields $\beta \in C_X(v)$. For the proof of (ii) note first that $C_X(v) \subset \bigcap_{f \in A \setminus \{0\}} \mathrm{Cont}^{\geq v(f)}(f)$. Let $C'$ be an irreducible component of $\bigcap_{f \in A \setminus \{0\}} \mathrm{Cont}^{\geq v(f)}(f)$ containing $C_X(v)$. Let $\alpha$ and $\alpha'$ be the generic points of $C_X(v)$ and $C'$, respectively. As $\alpha' \in \bigcap_{f \in A \setminus \{0\}} \mathrm{Cont}^{\geq v(f)}(f)$, we have $\mathrm{ord}_t\, \alpha'^*(f) \geq v(f)$ for every $f \in A \setminus \{0\}$. On the other hand, the inclusion $\alpha \in C'$ and the upper-semicontinuity yield $\mathrm{ord}_t\, \alpha'^*(f) \leq \mathrm{ord}_t\, \alpha^*(f) = v(f)$ for every $f \in A \setminus \{0\}$. Therefore we obtain that $v_{\alpha'} = v$ which implies $\alpha' \in C_X(v)$. For the proof of (iii), let $\beta$ be the generic point of $C_X(v')$. Then, $v|_A \leq v'|_A$ if and only if $\mathrm{ord}_t\, \beta^*(f) = v'(f) \geq v(f)$ for every $f \in A \setminus \{0\}$ and this is equivalent to the inclusion $\beta \in \bigcap_{f \in A \setminus \{0\}} \mathrm{Cont}^{\geq v(f)}(f)$. $\square$

The proposition suggests that $v|_A \leq v'|_A$ does not imply $C_X(v) \supset C_X(v')$, since there is no guarantee of the equality $\overline{\bigcap_{f \in A \setminus \{0\}} \mathrm{Cont}^{v(f)}(f)} = \bigcap_{f \in A \setminus \{0\}} \mathrm{Cont}^{\geq v(f)}(f)$. In the following we show an example for which the equality actually does not hold.

In the rest of this section, $X = \mathbb{C}^2 = \operatorname{Spec} A$ and $A = \mathbb{C}[x,y]$. We construct an example of $v, v'$ such that $v|_A \leq v'|_A$ but $C_X(v) \not\supset C_X(v')$.

**4.2 (Construction).** Let $X_\infty = \operatorname{Spec} \mathbb{C}[a_i, b_i \mid i \geq 0]$. Then
$$\operatorname{Cont}^2(x) \cap \operatorname{Cont}^3(y) = \operatorname{Spec} \mathbb{C}[a_2^{\pm 1}, a_3, .., b_3^{\pm 1}, b_4, ..].$$
Let $D \subset \operatorname{Cont}^2(x) \cap \operatorname{Cont}^3(y)$ be the closed subset defined by the ideal $(a_2^3 - b_3^2, 3a_2^2 a_3 - 2b_3 b_4)$ in $\operatorname{Cont}^2(x) \cap \operatorname{Cont}^3(y)$. Then, it is easy to check that $D$ is irreducible. Let $\alpha \in D$ be the generic point and $v$ be the divisorial valuation $v_\alpha$ corresponding to $\alpha$. Let $v'$ be the toric valuation with $v'(x) = 3, v'(y) = 4$

**Theorem 4.3.** *Let $v, v'$ be as above. Then, the following hold:*
  (i) $\alpha^*(x^3 - y^2) = c_8 t^8 + c_9 t^9 + \ldots$, where $c_8$ is transcendental over $\mathbb{C}(a_2, a_3, b_3, b_4)$ which is a subfield of $K$, where $\alpha^* : A \longrightarrow K[[t]]$ is a ring homomorphism corresponding to $\alpha$;
  (ii) Let $Z \subset X$ be the closed subset defined by $x^3 - y^2 = 0$ and $\psi_m : X_\infty \longrightarrow X_m$ the canonical projection for $m \in \mathbb{N}$. Then,
$$D = \operatorname{Cont}^2(x) \cap \operatorname{Cont}^3(y) \cap \psi_7^{-1}(Z_7) = \operatorname{Cont}^2(x) \cap \operatorname{Cont}^3(y) \cap \operatorname{Cont}^{\geq 8}(x^3 - y^2);$$
  (iii) $C_X(v) = \overline{D}$;
  (iv) $v|_A \leq v'|_A$;
  (v) $C_X(v) \not\supset C_X(v')$.

*Proof.* (i) Note that $K$ is the quotient field of the ring $\mathbb{C}[a_2^{\pm 1}, a_3, \ldots, b_3^{\pm 1}, b_4, ..]/(a_2^3 - b_3^2, 3a_2^2 a_3 - 2b_3 b_4)$. The image $\alpha^*(x^3 - y^2)$ is
$$(a_2^3 - b_3^2)t^6 + (3a_2^2 a_3 - 2b_3 b_4)t^7 + (3a_2 a_3^2 + 3a_2^2 a_4 - b_4^2 - 2b_3 b_5)t^8 + \ldots,$$
where the coefficients of $t^6$ and $t^7$ are zero in $K$ and the coefficient of $t^8$ contains $a_4$ and $b_5$ which are algebraically independent over the subfield $\mathbb{C}(a_2, a_3, b_3, b_4)$ in $K$. Therefore $c_8$ is transcendental over $\mathbb{C}(a_2, a_3, b_3, b_4)$.

(ii) Let $f_j$ be the coefficient of $t^j$ in
$$(a_0 + a_1 t + ..)^3 - (b_0 + b_1 t + \ldots)^2.$$
Then, the closed subset $\psi_7^{-1}(Z_7)$ is defined by the equations $f_0(a_i, b_i) = \ldots = f_7(a_i, b_i) = 0$ in $X_\infty$. Noting that $a_0 = a_1 = b_0 = b_1 = b_2 = 0$ in $\operatorname{Cont}^2(x) \cap \operatorname{Cont}^3(y)$, we obtain that $f_j = 0$ automatically holds for $j = 0, .., 5$ in $\operatorname{Cont}^2(x) \cap \operatorname{Cont}^3(y)$. On the other hand $f_6$ and $f_7$ coincide with $a_2^3 - b_3^2$ and $3a_2^2 a_3 - 2b_3 b_4$, respectively in $\operatorname{Cont}^2(x) \cap \operatorname{Cont}^3(y)$. Therefore, $D = \operatorname{Cont}^2(x) \cap \operatorname{Cont}^3(y) \cap \psi_7^{-1}(Z_7)$.

(iii) As $v_\alpha = v$, it is clear that $\overline{D} \subset C_X(v)$. On the other hand, by Proposition 4.1,

$$C_X(v) = \overline{\bigcap_{f \in A} \operatorname{Cont}^{v(f)}(f)}$$

and the right hand side is contained in $\operatorname{Cont}^2(x) \cap \operatorname{Cont}^3(y) \cap \operatorname{Cont}^{\geq 8}(x^3 - y^2)$.

For (iv) and (v) we need the following lemmas.

**Lemma 4.4.** *Let $v_0$ be the toric valuation with $v_0(x) = 2$ and $v_0(y) = 3$. Let $g \in A$ be a homogeneous element with respect to $v_0$, i.e., every monomial in $g$ has the same value of $v_0$. Then, $g = (x^3 - y^2)^m g'$, where $m \geq 0$ and $v(g') = v_0(g')$. Here, $g$ satisfies $v(g) > v_0(g)$ if and only if $m > 0$.*

*Proof.* First, if $v(g) = v_0(g)$, then $g$ is written as in the lemma with $m = 0$. Now assume that $v(g) > v_0(g)$. Note that $\alpha^*(x) = a_2 t^2 + a_3 t^3 + \ldots$, and $\alpha^*(y) = b_3 t^3 + b_4 t^4 + \ldots$, where $\alpha^* : A \longrightarrow K[[t]]$ is the ring homomorphism corresponding to $\alpha$. Here, $K$ is the quotient field of $\mathbb{C}[a_2^{\pm 1}, a_3, .., b_3^{\pm 1}, b_4, ..]/(a_2^3 - b_3^2, 3a_2^2 a_3 - 2b_3 b_4)$. As $g$ is homogeneous with respect to $v_0$, we obtain that

$$\alpha^*(g) = g(a_2, b_3) t^{v_0(g)} + (\text{higher order terms}).$$

Since $\operatorname{ord}_t \alpha^*(g) = v(g) > v_0(g)$, the leading coefficient $g(a_2, b_3) = 0$ in $K$. Therefore, $g(a_2, b_3) = 0$ in $\mathbb{C}[a_2, b_3]/(a_2^3 - b_3^2) \subset K$, which yields that $(x^3 - y^2) \mid g$. Now write $g = (x^3 - y^2) g'$. If $g'$ satisfies still $v(g') > v_0(g')$, then, by the above discussion, it follows that $(x^3 - y^2) \mid g'$. We have $g = (x^3 - y^2)^2 g''$ and $v_0(g'') < v_0(g') < v_0(g)$. By this procedure, we obtain finally the statement of the lemma. □

**Lemma 4.5.** *Assume $v(g) > v_0(g)$ for a polynomial $g \in A$. Let $g = g_1 + g_2 + \ldots + g_r$ be the decomposition into homogeneous parts with $v_0(g_i) < v_0(g_{i+1})$. Then,*

$$v(g) = \min\{v(g_i) \mid i = 1, \ldots, r\}.$$

*Proof.* By Lemma 4.4 we obtain for each $i = 1, .., r$

$$g_i = (x^3 - y^2)^{n_i} g_i',$$

where $n_i \geq 0$ and $g_i'$ satisfies $v(g_i') = v_0(g_i')$. As $v(x^3 - y^2) = 8$, it follows that $v(g_i) = 8n_i + v_0(g_i')$. Let $d = \min\{v(g_i) \mid i = 1, .., r\}$ and $I = \{i \mid v(g_i) = d\}$. Then, for two distinct element $i, j \in I$, it follows that $n_i \neq n_j$. Indeed, if $n_i = n_j$, then, $8n_i + v_0(g_i') = d = 8n_j + v_0(g_j')$, and therefore $v_0(g_i') = v_0(g_j')$. Hence, $v_0(g_i) = 6n_i + v_0(g_i') = 6n_j +$

$v_0(g'_j) = v_0(g_j)$, which is a contradiction to the assumption. Now, we have

$$\alpha^*(g) = \sum_{i \in I} \alpha^* \left( (x^3 - y^2)^{n_i} g'_i \right) + \sum_{i \notin I} \alpha^* \left( (x^3 - y^2)^{n_i} g'_i \right).$$

The right hand side is

$$\left( \sum_{i \in I} c_8^{n_i} g'_i(a_2, b_3) \right) t^d + \text{(higher order term)}.$$

If $v(g) = \mathrm{ord}_t \alpha^*(g) > d$, then $\#I \geq 2$ and we have a non-trivial algebraic relation

$$\sum_{i \in I} c_8^{n_i} g'_i(a_2, b_3) = 0$$

of $c_8$ over $\mathbb{C}(a_2, b_3)$. But this is a contradiction to (i). Therefore, $v(g) = d$. □

*Proof of (iv) and (v) of Theorem 4.3.* Take an arbitrary element $g \in A$. Let $g = g_1 + .. + g_r$ be the homogeneous decomposition as in Lemma 4.5. Let $g_i = (x^3 - y^2)^{n_i} g'_i$, where $v(g'_i) = v_0(g'_i)$ as in Lemma 4.4. Then $v'(g_i) \geq 8n_i + v_0(g'_i) = v(g_i)$. By Lemma 4.5, $v(g) = \min\{v(g_i) \mid i = 1, .., r\} \leq \min\{v'(g_i) \mid i = 1, .., r\} \leq v'(g)$. This proves (iv).

Let $\beta \in \mathrm{Cont}^3(x) \cap \mathrm{Cont}^4(y)$ be the generic point. Then $\beta^*(x) = a_3 t^3 + a_4 t^4 + \ldots$ and $\beta^*(y) = b_4 t^4 + b_5 t^5 + \ldots$ Then

$$\beta \in D' = \mathrm{Spec}\, \mathbb{C}[a_2, a_3, .., b_3, b_4 \ldots]/(a_2^3 - b_3^2, 3a_2^2 a_3 - 2b_3 b_4) \subset X_\infty.$$

Here, $D'$ has two component. One is $\overline{D}$ and we denote the other by $D''$. The component $D''$ is defined by $a_2 = b_3 = 0$ in $D'$. It is easy to see that $\beta$ is the generic point of $D''$. Therefore, $C_X(v) = \overline{D}$ and $C_X(v') = D''$ have no inclusion relation. □

**4.6.** This theorem and 3.11, (i) shows that the relation $C_X(v) \supset C_X(v')$ is strictly stronger than the relation $v|_A \leq v'|_A$. Now we consider another relation of divisorial valuation. For $X = \mathbb{C}^2$ we have an order $\prec$ of divisors $E$ and $E'$ over $X$: $E \prec E'$ if there is a successive blowing-ups $\longrightarrow X_j \longrightarrow \cdots \longrightarrow X_i \longrightarrow X_{i-1} \longrightarrow \cdots X_1 \longrightarrow X$ such that $E$ appears on $X_i$ and $E'$ appears on $X_j$ with $i < j$, and the center of $E'$ on $X_i$ is contained in $E$.

**Proposition 4.7.** *Let $v = \mathrm{val}_E$ and $v' = \mathrm{val}_{E'}$. If $E \prec E'$, then $C_X(v) \supset C_X(v')$. But $C_X(v) \supset C_X(v')$ does not imply $E \prec E'$.*

*Proof.* Let $p \in E \subset X_i$ be the center of $E'$ on $X_i$. Take a suitable affine neighborhood $U$ of $p$, then $U = \mathrm{Spec}\, B \simeq \mathbb{C}^2$, $p$ is the origin and $E$ is an invariant divisor, where we put a suitable toric structure on $U$. Then we have $v|_B \leq v'|_B$. As $v$ is toric on $U$, this inequality implies $C_U(v) \supset C_U(v')$ by Lemma 3.11 and therefore $C_X(v) \supset C_X(v')$ by Proposition 2.9.

For the second assertion, let $v = (2, 1)$ and $v' = (2, 3)$ in $\sigma \cap N$. We also denote the toric valuations corresponding to $v$ and $v'$ by the same symbol $v$ and $v'$, respectively. Then $v|_A \leq v'|_A$, where $A = \mathbb{C}[\sigma^\vee \cap M]$. Since $v$ and $v'$ are toric, it follows that $C_X(v) \supset C_X(v')$. But the divisors $E$ and $E'$ corresponding to $v$ and $v'$, respectively, does not satisfy $E \prec E'$. Indeed, let $\varphi_1 : X_1 \longrightarrow X$ be the blow-up at 0 and $E_1$ the exceptional divisor. There are two closed orbits $p_1$ and $p_2$ on $E_1$. The divisor $E$ is the exceptional divisor of the blow-up $\varphi : X_2 \longrightarrow X_1$ at one of the closed orbits on $E_1$, let it be $p_1$. Then, the center of $E'$ on $X_2$ is $p_2$. Therefore $E \not\prec E'$. $\square$

By this, the relation of $v$ and $v'$ for $C_X(v) \supset C_X(v')$ is something between "$v \prec v'$" and "$v|_A \leq v'|_A$".

Department of Mathematics, Tokyo Institute of Technology, Oh-Okayama, Meguro, Tokyo, Japan
e-mail : shihoko@math.titech.ac.jp